\documentclass[12pt,psfig]{article}
\usepackage{graphicx,epsfig}
\usepackage{amsmath}
\usepackage{mathrsfs}
\usepackage{amssymb}
\def\ge{\geq}

\def\bbb{\begin{eqnarray*}}

\def\eee{\end{eqnarray*}}

\begin{document}

\baselineskip=17pt
\begin{center}

\vspace{-0.6in} {\large \bf Estimations of topological entropy for\\
non-autonomous discrete systems }
\\ [0.2in]

Hua Shao, Yuming Shi$^{*}$, Hao Zhu

\vspace{0.15in} Department of Mathematics, Shandong University \\
 Jinan, Shandong 250100, P.~R. China\\

\footnote{$^{*}$ The corresponding author.}
\footnote{$^{**}$ Email addresses: huashaosdu@163.com (H. Shao), ymshi@sdu.edu.cn (Y. Shi), haozhusdu@163.com (H. Zhu).} \

\end{center}

{\bf Abstract.} In this paper, an estimation of lower bound of topological entropy for coupled-expanding systems associated with transition matrices in compact Hausdorff spaces is given. Estimations of upper and lower bounds of topological entropy for systems in compact metric spaces are obtained by their topological equi-semiconjugacy to subshifts of finite type under certain conditions.\medskip

\medskip

{2010 {\bf \it Mathematics Subject Classification}}: 37B40, 37B55, 37B10.

\bigskip

\noindent{\bf 1. Introduction}\medskip

The concept of topological entropy was first introduced by Adler, Konhelm and McAndrew in the 1960s,
based on open covers for continuous maps in compact topological spaces [1].
In 1970, Bowen gave another definition based on separated and spanning sets for uniformly continuous
maps in metric spaces [5], and this definition is equivalent to Adler's definition for
continuous maps in compact metric spaces.

Topological entropy provides a numerical measure for the complexity of dynamical systems, and it has close
relationships with many important dynamical properties, such as chaos, Lyapunove exponents and so on.
Thus, a lot of attention has been focused on computations and estimations of topological
entropy of a dynamical system. Many good results have been obtained for autonomous discrete systems [2-5, 7, 8, 11, 15-17].

In 1996, Kolyada and Snoha extended the concept of topological entropy for autonomous discrete systems to non-autonomous discrete systems
and obtained a series of important properties of it [13]. They gave two definitions also based on open covers, and separated and spanning sets,
separately, and proved that these two definitions are equivalent for continuous maps in
compact metric spaces, too. They studied some relations of topological entropy between two topologically
equi-semiconjugate non-autonomous systems and proved that topological entropy is an invariant value of
topological equi-conjugacy. Their works were followed by some other scholars [6, 14, 25, 27], in which Zhu and his coauthors studied
relations between topological entropy and measure-theoretic entropy [27].

The theory of symbolic dynamical systems played a significant role in the study of dynamical systems.
In the 1960s, complex dynamical behaviors of Smale horseshoe were depicted by its topological conjugacy to the full
shift [24]. Later, many scholars studied complexity of dynamical systems in this way [3, 12, 19-23, 26].
Motivated by these results, we shall try to give estimations of topological entropy for non-autonomous discrete systems
by their topological equi-semiconjugacy to subshifts of finite type in the present paper.

The rest of the paper is organized as follows. In Section 2, some basic concepts and lemmas are presented.
In Sections 3 and 4, estimations of topological entropy for non-autonomous discrete systems are discussed.

\bigskip

\noindent{\bf 2. Preliminaries}\medskip

In this section, some basic concepts and useful lemmas are presented.\medskip

We shall consider the following non-autonomous discrete system in the present paper:
\vspace{-0.2cm}$$x_{n+1}=f_n(x_n),\;\;n\geq1,                             \ \ \ \ \eqno(2.1)\vspace{-0.2cm}$$
where $f_n: X\to X$ is a continuous map for each $n\ge 1$, and $X$ is a topology space.

Firstly, we shall briefly recall the definition of topological entropy of system (2.1)
based on open covers given in [13]. Assume that $X$ is compact.
For open covers $\mathscr{A}_{1},\cdots,\mathscr{A}_{n},\; \mathscr{A}$ of $X$, denote
$\bigvee_{i=1}^{n}\mathscr{A}_{i}:=\{\bigcap_{i=1}^{n}A_{i}: A_{i}\in\mathscr{A}_{i},\  1\leq i\leq n\}$,
$f_{i}^{-n}(\mathscr{A}):=\{f_{i}^{-n}(A): A\in\mathscr{A}\}$, and $\mathscr{A}_{i}^{n}:=\bigvee_{j=0}^{n-1}f_{i}^{-j}(\mathscr{A})$, where
$f_{i}^{n}=f_{i+n-1}\circ\cdots\circ f_{i+1}\circ f_{i}$, and $f_{i}^{-n}=f_{i}^{-1}\circ f_{i+1}^{-1}\circ\cdots\circ f_{i+n-1}^{-1}$. Let $\mathcal{N}(\mathscr{A})$ be the minimal possible cardinality of all subcovers chosen from $\mathscr{A}$.
Then the topological entropy of $f_{1,\infty}=\{f_{i}\}_{i=1}^{\infty}$ or system (2.1) on the cover $\mathscr{A}$ is defined by
\vspace{-0.2cm}$$h(f_{1, \infty},\mathscr{A}):=\limsup_{n\rightarrow\infty}\log\mathcal{N}(\mathscr{A}_{1}^{n})/n,\vspace{-0.2cm}$$
and the topological entropy of $f_{1,\infty}$ or system (2.1) is defined by
\vspace{-0.2cm}$$h(f_{1, \infty}):=\sup\{h(f_{1, \infty},\mathscr{A}): \mathscr{A}\; {\rm is\; an\; open\; cover\; of\; X\;}\}.\vspace{-0.2cm}$$

Note that the two definitions of topological entropy of system (2.1) based on separated sets and
open covers are equivalent in the case that $X$ is a compact metric space [13].

Let $Y$ be any nonempty subset of $X$.
Denote the cover $\{A\cap Y: A\in\mathscr{A}\}$ of the set $Y$ by $\mathscr{A}|_{Y}$.
Then the topological entropy of $f_{1,\infty}$ on the set $Y$ is defined by
\vspace{-0.2cm}$$h(f_{1,\infty},Y):=\sup\{\limsup_{n\rightarrow\infty}\log\mathcal{N}(\mathscr{A}_{1}^{n}|_{Y})/n: \mathscr{A}\; {\rm is\; an\; open\; cover\; of\; X\;}\}.\vspace{-0.2cm}$$

Let $\Lambda_n\subset X$, $n\geq 1$, be nonempty subsets. The sequence $\{\Lambda_n\}_{n=1}^{\infty}$
of sets is said to be invariant under system (2.1) if $f_n(\Lambda_n)\subset\Lambda_{n+1}$ for all $n\geq 1$ [20].
Then system (2.1) restricted to $\{\Lambda_n\}_{n=1}^{\infty}$ is said to be an invariant subsystem of it. It is also called to be
the invariant subsystem of system (2.1) on $\{\Lambda_n\}_{n=1}^{\infty}$. 

Next, we shall recall topologically equi-conjugacy (equi-semiconjugacy) between two no-autonomous systems. 
We consider another system:
\vspace{-0.2cm}$$u_{n+1}=g_n(u_n),\;\;n\geq1,                                \ \ \ \ \eqno(2.2)\vspace{-0.2cm}$$
where $g_n:\, Y\to Y$ is a continuous map for each $n\ge 1$, and $Y$ is a topology space. \medskip

\noindent{\bf Definition 2.1} [20, Definition 3.3]. Let $(X, d)$ and $(Y, \rho)$ be  metric spaces, and $\{\Lambda_n\}_{n=1}^{\infty}$ and $\{E_n\}_{n=1}^{\infty}$ be invariant under systems (2.1) and (2.2), respectively. If, for each $n\geq1$, there exists an equi-continuous surjective map $h_n: \Lambda_n\to E_n$
such that $h_{n+1}\circ f_n=g_n\circ h_{n}$, then the invariant subsystem of system (2.1) on $\{\Lambda_n\}_{n=1}^{\infty}$ is said to be topologically $\{h_n\}_{n=1}^{\infty}$-equi-semiconjugate to the invariant subsystem of system (2.2) on $\{E_n\}_{n=1}^{\infty}$. Furthermore, if $\{h_n^{-1}\}_{n=1}^{\infty}$ is also equi-continuous, they are said to be topologically $\{h_n\}_{n=1}^{\infty}$-equi-conjugate.\medskip

By the method used in the proof of Theorem B in [13], one can get the following result:\medskip

\noindent{\bf Lemma 2.1.} {\it Let $(X, d)$ and $(Y, \rho)$ be compact metric spaces, and $\{\Lambda_n\}_{n=1}^{\infty}$ and $\{E_n\}_{n=1}^{\infty}$ be sequences of closed sets in $X$ and $Y$, respectively.  If an invariant subsystem of system {\rm (2.1)} on $\{\Lambda_n\}_{n=1}^{\infty}$ is topologically equi-semiconjugate to an invariant subsystem of system {\rm (2.2)} on $\{E_n\}_{n=1}^{\infty}$, then $h(g_{1,\infty}, E_1)\leq h(f_{1,\infty},\Lambda_1)$.}\medskip

Now, we shall recall some properties of transition matrices and subshifts of finite type.
A matrix $A=(a_{ij})_{N\times N}$ is said to be a transition matrix
if $a_{ij}=0$ or 1 for all $i,j$; $\sum_{j=1}^{N}a_{ij}\geq 1$ for all $i$; and
$\sum_{i=1}^{N}a_{ij}\geq 1$ for all $j$, $1\leq i,j\leq N$ [9]. A finite
sequence $w=(s_1,\cdots,s_l)$ is called an allowable word of length $l+1$ for $A$
if $a_{s_{i}s_{i+1}}=1, \;1\leq i\leq l-1$, where $1\leq s_{1},\cdots,s_{l}\leq N$.
For any $k\geq 1$ and any $1\leq i, j \leq N$, there are exactly $a_{ij}^{(k)}$ allowable
words of length $k +1$ for $A$, starting at $i$ and ending at $j$, in the form of
$(i,s_{1},\cdots,s_{k-1},j)$, where $a_{ij}^{(k)}$ denotes the $(i,j)$ entry of $A^{k}$ [18].
For a given transition matrix $A=(a_{ij})_{N\times N}$, denote
$\Sigma_{N}^{+}(A):=\{s=(s_{0}, s_{1}, \cdots): 1\leq s_{j}\leq N,\ a_{s_{j}s_{j+1}}=1,\;j\geq0\}$.
The map $\sigma_{A}:\Sigma_{N}^{+}(A)\rightarrow \Sigma_{N}^{+}(A)$ with
$\sigma_A((s_{0}, s_{1},\cdots)):=(s_{1}, s_{2},\cdots)$ is said to 
be the subshift of finite type for matrix $A$. Its topological entropy is equal to $\log\rho(A)$, where $\rho(A)$ is the spectral radius of matrix $A$. It is known that 
$\rho(A)=\lim_{n\rightarrow\infty}\Vert A^{n}\Vert^{\frac{1}{n}}$, where $\Vert A\Vert=\sum_{1\leq i,j\leq N}a_{ij}$ [???].\medskip

\noindent{\bf Definition 2.2} [20, Definition 2.10]. Let $(X, d)$ be a metric space. Assume that $A=(a_{ij})_{N\times N}$ is a transition matrix ($N\geq2$). If there exist $N$ subsets $V_{i}$ of $X$ with
$V_{i}\cap V_{j}=\partial V_{i}\cap\partial V_{j}$ for each pair $(i, j) ,\; 1\leq i\neq j\leq N$, such that
$f_{n}(V_{i})\supset {\bigcup_{a_{ij}=1}}V_{j},\;n\geq 1, $
where $\partial V_{i}$ is the boundary of $V_{i}$ in $X$, then system (2.1)
is said to be $A$-coupled-expanding in $V_{i}$, $1\leq i\leq N$.
Further, system (2.1) is said to be strictly $A$-coupled-expanding in $V_{i}$, $1\leq i\leq N$, if $d(V_{i}, V_{j})>0$ for all $1\leq i\neq j\leq N$.
In the special case that $a_{ij}=1$ for all $1\leq i,j\leq N$,
it is briefly said to be coupled-expanding or strictly coupled-expanding in $V_{i},\; 1\leq i\leq N$.\medskip

Note that it is not necessarily required that $X$ is a metric space in the definition of $A$-coupled-expansion.
The following result will be
also used in the sequent sections.\medskip

\noindent{\bf Lemma 2.2}. {\it Let $(X, d)$ and $(Y_n, \rho)$ be metric spaces and $\pi_{n}: X\rightarrow Y_n$, $n\geq1$. If $X$ is compact and $\{\pi_n\}_{n=1}^{\infty}$ is equi-continuous at any point of $X$; that is, for any fixed point $x\in X$, for any $\epsilon>0$, there exists $\delta>0$ such that for any $y\in B(x,\delta)$, $\rho(\pi_{n}(x), \pi_{n}(y))<\epsilon$ for any $n\geq1$, then $\{\pi_n\}_{n=1}^{\infty}$ is equi-continuous in $X$.}\medskip

\noindent{\bf Proof.} Fix $\epsilon>0$. By the assumption, for any point $x\in X$, there exists $\delta_{x}>0$ such that for any $y\in B(x,\delta_{x})$, $\rho(\pi_{n}(x), \pi_{n}(y))<\epsilon$ for any $n\geq1$. Then $\{B(x, \delta_{x}/2): x\in X\}$ is an open cover of $X$. Suppose that $\{B(x_i,\delta_{x_i}/2): 1\leq i\leq m\}$ is one of its finite subcovers since $X$ is compact. Set $\delta:=\min_{1\leq i\leq m}\{\delta_{x_i}\}$. For any $y_{1}, y_2\in X$ with $d(y_1,y_2)<\delta/2$, there exists $1\leq i\leq m$ such that $y_1\in B(x_i,\delta_{x_i}/2)$, then
$d(y_2, x_i)\leq d(y_2, y_1)+d(y_1, x_i)<\delta/2+\delta_{x_{i}}/2\leq\delta_{x_{i}}$. So, $y_1, y_2\in B(x_i,\delta_{x_i})$. Thus, $\rho(\pi_{n}(y_1), \pi_{n}(y_2))<2\epsilon$ for any $n\geq1$. Hence, $\{\pi_n\}_{n=1}^{\infty}$ is equi-continuous in $X$.

\bigskip

\noindent{\bf 3. An estimation of topological entropy for $A$-coupled-expanding systems}\medskip

In this section, we shall estimate lower bounds of topological entropy for $A$-coupled-expanding systems.\medskip

\noindent{\bf Theorem 3.1.} {\it Let $X$ be a compact Hausdorff space. Assume that there exist $N$ disjoint closed subsets $V_{1},\cdots,V_{N}$ of $X$ and an $N\times N$ transition matrix $A$ such that system {\rm(2.1)} is $A$-coupled-expanding in $V_{i}, \;1\leq i\leq N$.
Then $h(f_{1,\infty})\geq \log\rho(A)$, where $\rho(A)$ is the spectral radius of $A$.}\medskip

\noindent{\bf Proof.} Since $X$ is a compact Hausdorff space and $V_{i}$, $1\leq i\leq N$, are disjoint and closed subsets of $X$,
there exist pairwise disjoint open subsets $G_{1},\cdots, G_{N}$ such that $V_{i}\subset G_{i}, \;1\leq i\leq N$.
Set $G_{N+1}:=X\setminus (\bigcup_{i=1}^{N}V_{i})$.
Then $\mathscr{A}:=\{G_{1},\cdots,G_{N+1}\}$ is an open cover of $X$.
Fix $n\geq1$. Denote $\Omega:=\{w: w$ is an allowable word of length $n$ for matrix $A$\}.
For any $w=(i_{1},\cdots, i_{n})\in\Omega$, set
\vspace{-0.2cm}$$V_{w}:=\bigcap_{k=1}^{n}f_{1}^{-k+1}(V_{i_{k}}),\; G_{w}:
=\bigcap_{k=1}^{n}f_1^{-k+1}(G_{i_k}).                                                 \eqno(3.1)\vspace{-0.2cm}$$
Then $V_{w}\neq\emptyset$ by the assumption that system (2.1) is $A$-coupled-expanding
in $V_{i}$, $1\leq i\leq N$.
Fix one point $x_{w}\in V_{w}$. Set
\vspace{-0.2cm}$$S:=\{x_{w}: w\in\Omega\}.                                    \ \ \ \ \eqno(3.2)\vspace{-0.2cm}$$

For any two allowable words $w=(i_{1},\cdots, i_{n}), w'=(i'_{1},\cdots ,i'_{n})\in\Omega$ with $w\neq w'$, namely, there exists $1\leq l\leq n$ such that $i_{l}\neq i'_{l}$, we shall show the following three assertions hold.

(i) $x_{w}\neq x_{w'}$.

If not, then $x_{w}=x_{w'}$, which implies that $f_{1}^{l-1}(x_{w})\in V_{i_{l}}\cap V_{i'_{l}}$ by (3.1). This is a contradiction because $V_{i_{l}}$ and $V_{i'_{l}}$ are disjoint.

(ii) $G_{w}$ and $G_{w'}$ are the unique open sets of $\mathscr{A}_{1}^{n}$ which includes $x_{w}$ and $x_{w'}$, respectively.

It is only need to show that the assertion is true for $G_w$.
Clearly, $x_{w}\in G_{w}$ by (3.1). Now, it is to show the uniqueness. Otherwise, there exists another subset $\bigcap_{k=1}^{n}f_1^{-k+1}(G_{j_k})$ such that $x_{w}\in\bigcap_{k=1}^{n}f_1^{-k+1}(G_{j_k})$. Then there exists $1\leq m\leq n$ such that $i_{m}\neq j_{m}$. We divide the rest of proof into the following two cases:

{\bf Case 1.} Suppose that $1\leq j_{m}\leq N$. Then $f_{1}^{m-1}(x_{w})\in G_{i_{m}}\cap G_{j_{m}}$, which is a contradiction since $G_{i_{m}}$ and
$G_{j_{m}}$ are disjoint.

{\bf Case 2.} Suppose that $j_{m}=N+1$. Then $1\leq i_{m}\leq N$ and $f_{1}^{m-1}(x_{w})\in V_{i_{m}}\cap G_{N+1}$ by (3.1),
which is a contradiction since $G_{N+1}\cap(\bigcup_{i=1}^{N}V_{i})=\emptyset$.

(iii) $G_{w}\cap G_{w'}=\emptyset$.

Otherwise, suppose that there exists $x\in G_{w}\cap G_{w'}$.
Then $f_{1}^{l-1}(x)\in G_{i_{l}}\cap G_{i'_{l}}$,
which is a contradiction since $G_{i_{l}}$ and $G_{i'_{l}}$ are disjoint.

Hence, by the above assertions (i)-(iii) one gets that the number of open sets in any subcover of $\mathscr{A}_{1}^{n}$ is not less than the number of members in $S$; that is, $\mathcal{N}(\mathscr{A}_{1}^{n})\geq |S|$, where $|S|$ denotes the cardinality of $S$. It follows from (i) and (3.2) that
$|S|=|\Omega|=\sum_{1\leq i,j\leq N}a_{ij}^{(n-1)}=\Vert A^{n-1}\Vert$.
So $\mathcal{N}(\mathscr{A}_{1}^{n})\geq\Vert A^{n-1}\Vert$. Thus one has that
\vspace{-0.2cm}$$h(f_{1,\infty}, \mathscr{A})=\limsup_{n\rightarrow\infty}\log \mathcal{N}(\mathscr{A}_{1}^{n})/n
\geq \limsup_{n\rightarrow\infty}\log\Vert A^{n-1}\Vert/n
=\log\rho(A).\vspace{-0.2cm}$$
Then $h(f_{1,\infty})\geq \log\rho(A)$.
This completes the proof. \medskip

\noindent{\bf Remark 3.1.} This theorem extends the result of Block with his coauthors
for continuous interval maps in [2] to non-autonomous discrete systems.\medskip

Note that $\rho(A)$ is not easily computed in general. The following result gives
a simpler estimation, which is a direct consequence of Theorem 3.1 and Theorem 8.1.22 in [9].

\medskip

\noindent{\bf Corollary 3.1.} {\it Let all the assumptions in Theorem 3.1 hold. Then
$h(f_{1,\infty})\geq\log\nu$, where $\nu=\max\{\min\limits_{i}\sum\limits_{j=1}^{N}a_{ij},
\min\limits_{j}\sum\limits_{i=1}^{N}a_{ij}\}$.
}\medskip

\bigskip

\noindent{\bf 4. Estimations of topological entropy for an invariant subsystem }\medskip

In this section, we shall estimate the topological entropy for an invariant subsystem.\medskip

\noindent{\bf Theorem 4.1.} {\it Let $(X,d)$ be a compact metric space. Assume that there exist $N$ closed subsets $V_{1},\cdots,V_{N}$ of $X$ and an $N\times N$ transition matrix $A=(a_{ij})_{N\times N}$ such that system {\rm(2.1)} satisfies the following conditions:
\begin{itemize}\vspace{-0.2cm}
\item [{\rm (i)}] $f_{n}(V_{i})\supset {\bigcup\limits_{a_{ij}=1}}V_{j},\;1\leq i, j\leq N,\; n\geq 1$\vspace{-0.2cm};
\item [{\rm (ii)}] for any given $\alpha=(a_{0},a_{1},\cdots)\in \Sigma _N^+(A)$, $d(V_{\alpha}^{m,n})$                                                                                uniformly converges to 0 with respect to $n\geq1$ as $m\rightarrow\infty$, where $V_{\alpha}^{m,n}:=\bigcap\limits_{k=0}^{m}f_{n}^{-k}(V_{a_{k}})$.\vspace{-0.2cm}
\end{itemize}
Then, for each $n\geq1$, there exists a nonempty compact set $\Lambda_{n}\subset X$ with $f_{n}(\Lambda_{n})=\Lambda_{n+1}$ and an equi-continuous surjective map $\pi_{n}:\Sigma_{N}^{+}(A)\to\Lambda_{n}$ such that $(\Sigma_{N}^{+}(A),\sigma_{A})$ is topologically $\{\pi_{n}\}_{n=1}^{\infty}$-equi-semiconjugate to
 the invariant subsystem of system {\rm (2.1)} on $\{\Lambda_n\}_{n=1}^{\infty}$. Consequently, $h(f_{1,\infty},\Lambda_{1})\leq\log\rho(A)$.}\medskip

Note that $V_{j}$, $1\leq j\leq N$, are not required to be disjoint here as that in Theorem 3.1.\medskip

\noindent{\bf Proof.} For each $n\geq1$, by assumption (i) and the continuity of $f_n$, it can be easily verified that for any $\alpha\in \Sigma_{N}^{+}(A)$ and $m\geq0$, $V_{\alpha}^{m,n}$ is a nonempty closed and bounded subset and satisfies that
$V_{\alpha}^{m,n}\supset V_{\alpha}^{m+1,n}$.
This, together with assumption (ii), yields that $\bigcap_{m=0}^{\infty}V_{\alpha}^{m,n}$
is a singleton set by Lemma 2.7 in [19]. Set
\vspace{-0.2cm}$$\{x^{n}(\alpha)\}:=\bigcap_{m=0}^{\infty}V_{\alpha}^{m,n},
\;\;\Lambda_{n}:=\{x^{n}(\alpha): \alpha\in\Sigma_{N}^{+}(A)\}.                          \eqno(4.1)\vspace{-0.2cm}$$
Clearly, $\Lambda_{n}\ne \emptyset $, $\Lambda_{n}\subset \bigcup _{i=1}^N V_i \subset X$,
and $f_{n}(x^{n}(\alpha))=x^{n+1}(\sigma_{A}(\alpha))$ for any $\alpha\in\Sigma_{N}^{+}(A)$,
$n\geq1$.
Hence, it follows from the fact that $\sigma_{A}$ is surjective that
$f_{n}(\Lambda_{n})=\Lambda_{n+1}$, $n\geq1$.

Define a map $\pi_{n}:\Sigma _N^+(A)\rightarrow\Lambda_{n}$ by $\pi_{n}(\alpha)=x^{n}(\alpha)$ for any $\alpha\in\Sigma_{N}^{+}(A)$. It is evident that $\pi_{n}$ is well defined and surjective.
In addition, we have that
\vspace{-0.2cm}$$f_{n}\circ \pi_{n}(\alpha)=f_{n}(x^{n}(\alpha))=x^{n+1}(\sigma_{A}(\alpha))=\pi_{n+1}\circ\sigma_{A}(\alpha),\; \alpha\in\Sigma _N^+(A), \;n\geq1.\vspace{-0.2cm}$$
Therefore, $\pi_{n+1}\circ\sigma_{A}=f_{n}\circ \pi_{n}$, $n\geq1$.

Fix any $\alpha=(a_{0},a_{1},\cdots)\in\Sigma _N^+(A)$. By assumption (ii), for any $\epsilon>0$, there exists $N_0>0$, such that
for any $m\geq N_0$, $d(V_{\alpha}^{m,n})<\epsilon$ for all $n\geq1$.
Set $\delta:=1/2^{N_0+1}$. For each $\beta=(b_{0},b_{1},\cdots)\in\Sigma _N^+(A)$ with $\rho(\alpha,\beta)<\delta$, then $a_{j}=b_{j},\;0\leq j\leq N_0+1$. So, $x^{n}(\alpha), x^{n}(\beta)\in V_{\alpha}^{N_0+1,n}$, and $d(\pi_{n}(\alpha),\pi_{n}(\beta))=d(x^{n}(\alpha),x^{n}(\beta))\leq d(V_{\alpha}^{N_0+1,n})<\epsilon$ for all $n\geq1$. Thus, $\{\pi_{n}\}_{n=1}^{\infty}$ is equi-continuous at $\alpha$. Hence, $\{\pi_{n}\}_{n=1}^{\infty}$ is equi-continuous in $\Sigma _N^+(A)$ by Lemma 2.2. Further, for any $n\geq1$, $\Lambda_{n}$ is compact since $\Sigma _N^+(A)$ is compact and $\Lambda_{n}=\pi_{n}(\Sigma _N^+(A))$.

Hence, $(\Sigma_{N}^{+}(A),\sigma_{A})$ is topologically $\{\pi_{n}\}_{n=1}^{\infty}$-equi-semiconjugate to the invariant subsystem of system {\rm (2.1)} on  $\{\Lambda_n\}_{n=1}^{\infty}$. By Lemma 2.1, $h(f_{1,\infty},\Lambda_{1})\leq h(\sigma_{A})=\log\rho(A)$. The proof is complete.\medskip

Assumption (ii) in the above theorem is not easily checked in general. Therefore, we give a verifiable condition as
follows.\medskip

\noindent{\bf Theorem 4.2.} {\it Let all the assumptions in Theorem 4.1 hold except that assumption
{\rm(ii)} is replaced by
\begin{itemize}\vspace{-0.2cm}
\item [{\rm (ii$_{a}$)}] there exists a constant $\lambda>1$ such that
\vspace{-0.2cm}$$d(f_{n}(x),f_{n}(y))\geq\lambda d(x,y), \forall\  x,y\in V_{i},\;1\leq i\leq N,\;n\geq1.\vspace{-0.2cm}$$
\end{itemize}
Then, all the results in Theorem 4.1 hold.}\medskip

\noindent{\bf Proof.} By the method used in the proof of Theorem 4.2 in [23], one can obtain that the assumption (ii)$_{a}$ implies assumption (ii) in Theorem 4.1. Hence, all the results in Theorem 4.1 hold and this completes the proof.\medskip

The following result gives an estimation of upper bound of topological entropy for the full system (2.1).\medskip

\noindent{\bf Theorem 4.3.} {\it Let all the assumptions in Theorem 4.2 hold and $\bigcup_{i=1}^{N}V_{i}=X$ except that assumption {\rm(i)} is replaced by
\begin{itemize}\vspace{-0.2cm}
\item [{\rm (i$_{a}$)}] $f_{n}(V_{i})= \bigcup\limits_{a_{ij}=1}V_{j},\;\; 1\leq i\leq N,\;\; n\geq 1.$
\end{itemize}\vspace{-0.2cm}
Then, for each $n\geq1$, there exists an equi-continuous surjective map $\pi'_{n}:\Sigma_{N}^{+}(A)\to X$ such that $(\Sigma_{N}^{+}(A),\sigma_{A})$ is topologically $\{\pi'_{n}\}_{n=1}^{\infty}$-equi-semiconjugate to system {\rm(2.1)}. Thus, $h(f_{1,\infty})\leq\log\rho(A)$.}\medskip

\noindent{\bf Proof.} By the proof of Theorems 4.1 and 4.2, it suffices to prove that $\Lambda_{n}=X$, $n\geq1$, where $\Lambda_{n}$ is specified in (4.1). In fact, for each $n\geq1$, by assumption (i$_{a}$) and $\bigcup_{j=1}^{N}V_{j}=X$ we have that for any $x\in X$, there exists $\beta=(b_{0}, b_{1},\cdots)\in\Sigma_{N}^{+}(A)$ such that $x=x^{n}(\beta)\in\Lambda_{n}$, which implies that $X\subset\Lambda_{n}$. Hence, $\Lambda_{n}=X$, $n\geq1$. The proof is complete.\medskip

\noindent{\bf Theorem 4.4.} {\it Let all the assumptions in Theorem 4.1 hold and $V_{1},\cdots,V_{N}$ are pairwise disjoint except that assumption
{\rm(ii)} is replaced by
\begin{itemize}\vspace{-0.2cm}
\item [{\rm (ii$_{b}$)}] there exists a positive constant $\mu$ such that
\vspace{-0.2cm}$$d(f_{n}(x),f_{n}(y))\leq \mu d(x,y), \forall\; x,y\in V_{i},\;1\leq i\leq N,\;n\geq1.\vspace{-0.2cm}$$
\end{itemize}
Then, for each $n\geq1$, there exists a nonempty compact subset $\Lambda_{n}\subset X$ with $f_{n}(\Lambda_{n})=\Lambda_{n+1}$ and an equi-continuous surjective map $h_{n}:\Lambda_{n}\to\Sigma_{N}^{+}(A)$ such that  the invariant subsystem of system {\rm (2.1)} on $\{\Lambda_n\}_{n=1}^{\infty}$ is topologically $\{h_{n}\}_{n=1}^{\infty}$-equi-semiconjugate to $(\Sigma_{N}^{+}(A),\sigma_{A})$. Consequently, $h(f_{1,\infty},\Lambda_{1})\geq\log\rho(A)$.}\medskip

\noindent{\bf Proof.} For each $n\geq1$, by the proof of Theorem 4.1, for any $\alpha\in \Sigma_{N}^{+}(A)$ and $m\geq0$, $V_{\alpha}^{m,n}$ is a nonempty closed subset and satisfies that $V_{\alpha}^{m,n}\supset V_{\alpha}^{m+1,n}$, where $V_{\alpha}^{m,n}$ is specified in Theorem 4.1. Then $\bigcap_{m=0}^{\infty}V_{\alpha}^{m,n}\neq\emptyset$ since $X$ is compact. Denote $\Lambda_{n}:=\bigcup_{\alpha\in\Sigma_{N}^{+}(A)}\bigcap_{m=0}^{\infty}V_{\alpha}^{m,n}$.
Then, $\Lambda_{n}$ is closed by the method used in the proof of Theorem 3.1 in [26], so is compact.

For any $x\in\Lambda_{n}$, there exists $\alpha\in\Sigma_{N}^{+}(A)$ such that $x\in\bigcap_{m=0}^{\infty}V_{\alpha}^{m,n}$. Define $h_{n}(x)=\alpha$, then the map $h_{n}: \Lambda_{n}\to\Sigma_{N}^{+}(A)$ is well defined and surjective since $V_{1},\cdots,V_{N}$ are disjoint. Further, it can be easily verified that
$f_{n}(x)\in \bigcap_{m=0}^{\infty}V_{\sigma_{A}(\alpha)}^{m,n+1}$.
So, $f_{n}(\Lambda_{n})\subset\Lambda_{n+1}$. On the other hand, by assumption (i) and the fact that $\sigma_{A}$ is surjective, one can obtain that $\Lambda_{n+1}\subset f_{n}(\Lambda_{n})$. Thus, $f_{n}(\Lambda_{n})=\Lambda_{n+1}$, $n\geq1$. In addition, one can also get that $h_{n+1}(f_n(x))=\sigma_{A}(\alpha)=\sigma_{A}(h_{n}(x))$. Therefore, $h_{n+1}\circ f_{n}=\sigma_{A}\circ h_{n}$.

With a similar argument to that used in the proof of Theorem 4.1 in [20], one can show that $\{h_n\}_{n=1}^{\infty}$ is equi-continuous.

Hence,  the invariant subsystem of system {\rm (2.1)} on $\{\Lambda_n\}_{n=1}^{\infty}$ is topologically $\{h_n\}_{n=1}^{\infty}$-equi-semiconjugate to $(\Sigma_{N}^{+}(A),\sigma_{A})$. Therefore, $h(f_{1,\infty},\Lambda_{1})\geq\log\rho(A)$ by Lemma 2.1. This completes the proof.\medskip

By Theorems 4.2 and 4.4, one can easily obtain the following result:\medskip

\noindent{\bf Theorem 4.5.} {\it Let all the assumptions in Theorem 4.4 and {\rm(ii$_{a}$)} hold, where $1<\lambda\leq\mu$.
Then, for each $n\geq1$, there exists a nonempty compact subset $\Lambda_{n}\subset X$ with $f_{n}(\Lambda_{n})=\Lambda_{n+1}$ and an equi-continuous homeomorphism $h_{n}:\Lambda_{n}\rightarrow\Sigma_{N}^{+}(A)$ such that  the invariant subsystem of system {\rm (2.1)} on $\{\Lambda_n\}_{n=1}^{\infty}$ is topologically $\{h_{n}\}_{n=1}^{\infty}$-equi-conjugate to $(\Sigma_{N}^{+}(A),\sigma_{A})$. Consequently, $h(f_{1,\infty},\Lambda_{1})=\log\rho(A)$.
}\medskip

\noindent{\bf Remark 4.1.} Let all the assumptions in Theorem 4.5 hold. If $A=(a_{ij})_{N\times N}$ is irreducible with
$\sum_{j=1}^{N}a_{i_0j}\geq2$ for some $1\leq i_0\leq N$, then $\rho(A)>1$, and consequently $h(f_{1,\infty},\Lambda_{1})>0$.
It was shown that system (2.1) is chaotic in the strong sense of Li-Yorke in this case by Theorem 4.1 in [20]. More recently, we relaxed the conditions of Theorem 4.1 in [20]; that is, we removed assumption (ii$_{b}$),
and got the same conclusion (see Theorem 4.2 in [23]).

\bigskip

\noindent{\bf \large References}
\def\hang{\hangindent\parindent}
\def\textindent#1{\indent\llap{#1\enspace}\ignorespaces}
\def\re{\par\hang\textindent}
\noindent \vskip 3mm

\re{[1]} R. Adler, A. Konheim, M. McAndrew. Toplogical entropy, Trans. Amer. Math. Soc. 114 (1965) 309--319.

\re{[2]} L. Block, J. Guckenheimer, M. Misuirewicz, L. S. Young. Periodic orbits and topological entropy of one-dimensional maps, Global Theory of Dynamical Systems, Lecture Notes in Math. vol. 819, Springer-Verlag, New York (1980) 18--34.

\re{[3]} L. Block, W. Coppel. Dynamics in One Dimension, Lecture Notes in Mathematics, 1513, Berlin/Heidelberg: Springer-Verlag, 1992.

\re{[4]} R. Bowen. Topological entropy and axiom A, Global Analysis, Proc. Sympos: Pure Math. 14 (1970) 23--41.

\re{[5]} R. Bowen. Entropy for group endomorphisms and homogenuous spaces, Trans. Amer. Math. Soc. 153 (1971) 401--414.

\re{[6]} J. S. C$\acute{a}$novas. On entropy of non-autonomous discrete systems, Prog. Chall. Dyn. Syst. 54 (2013) 143--159.

\re{[7]} L. W. Goodwyn. Topological entropy bounds and measure-theoretic Entropy, Proc. Am. Math. Soc. 23 (1969) 679--688.

\re{[8]} T. N. T. Goodman. Relating to entropy and measure entropy, Bull. Lond. Math. Soc. 3 ( 1971) 176--180.

\re{[9]} R. A. Horn, C. R. Johnson. Matrix Analysis. Vol. I, Cambridge Univ. Press, UK, 1985.

\re{[10]} Q. Huang, Y. Shi, L. Zhang. Sensitivity of non-autonomous discrete dynamical systems. Appl. Math. Lett. 39 (2015) 31--34.

\re{[11]} S. Ito. An Estimate from above for the Entropy and the topological entropy of a $C^{1}$ diffeomorphism, Proc. Japan  Acad. 46 (1970) 262--230.

\re{[12]} J. Kennedy, J. A. Yorke. Topological horseshoes, Trans. Amer. Math. Soc. 353 (2001) 2513--2530.

\re{[13]} S. Kolyada, L. Snoha. Topological entropy of non-autononous dynamical systems, Random Comp. Dyn. 4 (1996) 205--233.

\re{[14]} S. Kolyada, Mi. Misiurewicz, L. Snoha. Topological entropy of non-autonomous piecewise monotone dynamical systems on the interval,
  Fundam.  Math. 160 (1999) 161--181.

\re{[15]} A. Manning. Topological entropy and the first homology group, Dyn. Syst, Proc. Symp. Univ. Warwick 1973/74, Lect. Notes Math. 46 (1975) 185--190.

\re{[16]} M. Misiurewicz, F. Przytycki. Topological entropy and degree of smooth mappings, Bull. Acad. Pol. Sci. Math. Astron, Phys. 25 (1977) 573--574.

\re{[17]} C. Ri, H. Ju, X. Wu. Entropy for $A$-coupled-expanding maps and chaos. arXiv:1309.6769v2 [math.DS] 28 Sep (2013).

\re{[18]} C. Robinson. Dynamical Systems: Stability, Symbolic Dynamics and Chaos, Florida: CRC Press, Inc. 1999.

\re{[19]} Y. Shi, G. Chen. Chaos of discrete dynamical systems in complete metric spaces, Chaos Solit. Fract. 22 (2004) 555--571.

\re{[20]} Y. Shi, G. Chen. Chaos of time-varying discrete dynamical systems, J. Differ. Equ. Appl. 15 (2009) 429--449.

\re{[21]} Y. Shi, H. Ju, G. Chen. Coupled-expanding maps and one-sided symbolic dynamical systems, Chaos Solit. Fract. 39 (2009) 2138--2149.

\re{[22]} Y. Shi. Chaos in non-autonomous discrete dynamical systems approached by their sub-systems, Int. J. Bifurc. Chaos 22 (2012) 1250284.

\re{[23]} H. Shao, Y. Shi, H. Zhu. Strong Li-Yorke chaos for time-varying discrete systems with A-coupled-expansion, Submitted for Publication.

\re{[24]} S. Smale. Differentiable dynamical systems, Bull. Amer. Math. Soc. 73 (1967) 747--817.

\re{[25]} J. Zhang, L. Chen. Lower Bounds of the topological entropy of non-autonomous dynamical systems, Appl. Math. J. Chinese Univ, 24(1) (2009) 76--82.

\re{[26]} X. Zhang, Y. Shi, G. Chen, Some properties of coupled-expanding maps in compact sets, Proc. Amer. Math. Soc. 141(2) (2013) 585--595.

\re{[27]} Y. Zhu, Z. Liu, X. Xu,  W. Zhang. Entropy of non-autonomous dynamical systems, J. Korean Math. Soc. 49 (2012) 165--185.

\re{[28]} Z. Zhou. Symbolic Dynamics, Shanghai Scientific and Technological Education Publishing House, Shanghai, 1997.

\end{document}